\documentclass{article}
\usepackage{amssymb,amsmath,epsfig,amscd,eucal,psfrag}

\newtheorem{thm}{Theorem}[section]
\newtheorem{defn}[thm]{Definition}
\newtheorem{lem}[thm]{Lemma}
\newtheorem{prop}[thm]{Proposition}

\newtheorem{rem}[thm]{Remark}
\newtheorem{cor}[thm]{Corollary}
\newtheorem{exa}[thm]{Example}

\newenvironment{pf}{\par\medskip\noindent{\em Proof. }}{\hfill $\square$\par\medskip}
\newenvironment{pfof}[1]{\par\medskip\noindent{\em Proof of #1. }}{\hfill $\square$\par\medskip}

\newcommand{\F}{\mathbb{F}}

\newcommand{\Proj}{\mathbb{P}}
\newcommand{\R}{\mathbb{R}}
\newcommand{\Z}{\mathbb{Z}}

\newcommand{\ML}{\mathcal{ML}}

\newcommand{\curlyG}{\mathcal{G}}

\newcommand{\Axis}{\mathrm{Axis}}
\newcommand{\Hom}{\mathrm{Hom}}
\newcommand{\Aut}{\mathrm{Aut}}
\newcommand{\Mod}{\mathrm{Mod}}
\newcommand{\Stab}{\mathrm{Stab}}
\newcommand{\Isom}{\mathrm{Isom}}

\title{Solutions to Bestvina \& Feighn's Exercises on Limit Groups}
\date{16th March 2006}

\author{Henry Wilton}

\begin{document}

\maketitle

\begin{abstract}
This article gives solutions to the exercises in Bestvina and
Feighn's paper \cite{BF03} on Sela's work on limit groups.  We prove
that all constructible limit groups are limit groups and give an
account of the shortening argument of Rips and Sela.
\end{abstract}

Mladen Bestvina and Mark Feighn's beautiful first set of notes
\cite{BF03} on Zlil Sela's work on the Tarski problems (see
\cite{Se1} \emph{et seq.}) provides a very useful introduction to
the subject.  It gives a clear description of the construction of
Makanin--Razborov diagrams, and precisely codifies the structure
theory for limit groups in terms of \emph{constructible limit groups
(CLGs)}.  Furthermore, the reader is given a practical initiation in
the subject with exercises that illustrate the key arguments.  This
article is intended as a supplement to \cite{BF03}, to provide
solutions to these exercises.  Although we do give some definitions
in order not to interrupt the flow, we refer the reader to
\cite{BF03} for all the longer definitions and background ideas and
references.

\section{Definitions and elementary properties} \label{Introduction}

In this section we present solutions to exercises 2, 3, 4, 5, 6 and 7, which
give some of the simpler properties and the first examples and non-examples
of limit groups.

\subsection{$\omega$-residually free groups}
Fix $\F$ a free group of rank $r>1$.

\begin{defn}
A finitely generated group $G$ is \emph{$\omega$-residually free}
if, for any finite subset $X\subset G$, there exists a homomorphism
$h:G\rightarrow\F$ whose restriction to $X$ is injective.
(Equivalently, whenever $1\notin X$ there exists a homomorphism
$h:G\to\F$ so that $1\notin h(X)$.)
\end{defn}

Residually free groups inherit many of the properties of free
groups; the first and most obvious property is being torsion-free.

\begin{lem}[Exercise 2 of \cite{BF03}]\label{LGs are TF}
Any residually free group is torsion-free.
\end{lem}
\begin{pf}
Let $G$ be $\omega$-residually free (indeed $G$ can be thought of
as merely residually free).  Then for any $g\in G$, there exists a
homomorphism $h:G\rightarrow\F$ with $h(g)\neq 1$; so $h(g^k)\neq
1$ for all integers $k$, and $g^k\neq 1$.
\end{pf}

It is immediate that any subgroup of an $\omega$-residually free group is
$\omega$-residually free (exercise 6 of \cite{BF03}).

That the choice of $\F$ does not matter follows from the observation
that all finitely generated free groups are $\omega$-residually
free.

\begin{exa}[Free groups] \label{Free groups are LGs}
Let $F$ be a finitely generated free group.  Realize $\F$ as  the
fundamental group of a rose $\Gamma$ with $r$ petals; that is, the
wedge of $r$ circles.  Then $\Gamma$ has an infinite-sheeted cover
that corresponds to a subgroup $F'$ of $\F$ of countably infinite
rank. The group $F$ can be realized as a free factor of $F'$; this
exhibits an injection $F\hookrightarrow\F$.  In particular, every
free group is $\omega$-residually free.
\end{exa}

\begin{exa}[Free abelian groups] \label{Free abelian groups are LGs}
Let $A$ be a finitely generated free abelian group, and let
$a_1,\ldots,a_n\in A$ be non-trivial.  Fix a basis for $A$, and
consider the corresponding inner product.  Let $z\in A$ be such that
$\langle z,a_i\rangle\neq 0$ for all $i$.  Then inner product with
$z$ defines a homomorphism $A\to\Z$ so that the image of every $a_i$
is non-trivial, as required.
\end{exa}

Examples \ref{Free groups are LGs} and \ref{Free abelian groups
are LGs} give exercise 3 of \cite{BF03}.

\subsection{Limit groups}

Groups that are $\omega$-residually free are natural examples of
limit groups.

\begin{defn}
Let $\F$ be as above, and $\Gamma$ a finitely generated group.  A
sequence of homomorphisms $(f_n:\Gamma\rightarrow\F)$ is
\emph{stable} if, for every $g\in G$, $f_n(g)$ is either
eventually $1$ or eventually not $1$.  The \emph{stable kernel} of
a stable sequence of homomorphisms $(f_n)$ consists of all $g\in
G$ with $f_n(g)$ eventually trivial;  it is denoted
$\underrightarrow{ker}f_n$.

A \emph{limit group} is a group arising as a quotient
$\Gamma/\underrightarrow{ker}f_n$ for $(f_n)$ a stable sequence.
\end{defn}

\begin{lem}[Exercise 5 of \cite{BF03}]\label{ORFs are LGs}
Every $\omega$-residually free group is a limit group.
\end{lem}
\begin{pf}
Let $G$ be an $\omega$-residually free group.  Fix a generating
set, and let $X_n\subset G$ be the ball of radius $n$ about the
identity in the word metric.  Let $f_n:G\rightarrow\F$ be a
homomorphism that is injective on $X_n$.  Now $f_n$ is a stable
sequence and the stable kernel is trivial, so $G$ is a limit
group.
\end{pf}

In fact every limit group is $\omega$-residually free (lemma 1.11 of
\cite{BF03}).  Henceforth, we shall use the terms interchangeably.

\subsection{Negative examples}

Let's see some examples of groups that aren't limit groups.   The
first three examples are surface groups that aren't even
residually free.  It follows from lemma \ref{LGs are TF} that the fundamental group of the real projective plane is not a limit group.  A slightly finer analysis yields some other negative examples.

\begin{lem}\label{2-generator LGs}
The only 2-generator residually free groups are the free group of
rank 2 and the free abelian group of rank 2.
\end{lem}
\begin{pf}
Let $G$ be a residually free group generated by $x$ and $y$.  If
$G$ is non-abelian then $[x,y]\neq 1$ so there exists a
homomorphism $f:G\rightarrow\F$ with $f([x,y])\neq 1$.  So $f(x)$
and $f(y)$ generate a rank 2 free subgroup of $\F$.  Therefore $G$
is free.
\end{pf}

In particular, the fundamental group of the Klein bottle is not a
limit group.   Only one other surface group fails to be
$\omega$-residually free. This was first shown by R.~S.~Lyndon in \cite{L59}.

\begin{lem}[The surface of Euler characteristic -1]\label{Lyndon's example}
Let $\Sigma$ be the closed surface of Euler characteristic -1.
Then any homomorphism $f_*:\pi_1(\Sigma)\to\F$ has abelian image.
In particular, since $\pi_1(\Sigma)$ is not abelian, it is not
residually free.
\end{lem}
\begin{pf}
Let $\Gamma$ be a bouquet of circles so $\F=\pi_1(\Gamma)$.
Realize the homomorphism $f_*$ as a map from $\Sigma$ to $\Gamma$,
which we denote by $f$.  Our first aim is to find an essential
simple closed curve in the kernel of $f_*$.

Consider $x$ the mid-point of an edge of $\Gamma$.  Altering $f$
by a homotopy, it can be assumed that $f$ is transverse at $x$; in
this case, $f^{-1}(x)$ is a collection of simple closed curves.
Let $\gamma$ be such a curve. If $\gamma$ is null-homotopic in
$\Sigma$ then $\gamma$ can be removed from $f^{-1}(x)$ by a
homotopy.  If all components of the pre-images of all midpoints
$x$ can be removed in this way then $f_*$ was the trivial
homomorphism.  Otherwise, any remaining such component $\gamma$
lies in the kernel of $f$ as required.

We proceed with a case-by-case analysis of the components of
$\Sigma\smallsetminus\gamma$.
\begin{enumerate}
\item If $\gamma$ is 2-sided and separating then, by examining Euler characteristic,
the components of $\Sigma\smallsetminus\gamma$ are a punctured torus
or a Klein bottle, together with a M\"{o}bius band.  So $f$ factors
through the one-point union $T\vee \R P^2$ or $K\vee\R P^2$.  In
either case, it follows that the image is abelian.
\item If $\gamma$ is 2-sided and non-separating then $\Sigma-\gamma$ is the
non-orientable surface with Euler characteristic -1 and two boundary
components, so $f_*$ factors through $\Z*\Z/2\Z$ and hence through
$\Z$.
\item If $\gamma$ is 1-sided then $\gamma^2$ is 2-sided and separating, and
case 1 applies.
\end{enumerate}
This finishes the proof.
\end{pf}

Lemmas \ref{LGs are TF}, \ref{2-generator LGs} and \ref{Lyndon's
example} give exercise 4 of \cite{BF03}. Here is a more interesting
obstruction to being a limit group. A group $G$ is \emph{commutative
transitive} if every non-trivial element has abelian centralizer;
equivalently, if $[x,y]=[y,z]=1$ then $[x,z]=1$. Note that $\F$ is
commutative transitive, since every non-trivial element has cyclic
centralizer.

\begin{lem}[Exercise 7 of \cite{BF03}] \label{Abelian centralizers}
Limit groups are commutative transitive.
\end{lem}
\begin{pf}
Let $G$ be $\omega$-residually free, let $g\in G$, and suppose
$a,b\in G$ commute with $g$.  Then there exists a homomorphism
$$
f:G\rightarrow\F
$$
injective on the set $\{1,g,[a,b]\}$.  Then $f([g,a])=f([g,b])=1$;
since $\F$ is commutative transitive, it follows that
$$
f([a,b])=1.
$$
So $[a,b]=1$, as required.
\end{pf}

A stronger property also holds.  A subgroup $H\subset G$ is
\emph{malnormal} if, whenever $g\notin H$, $gHg^{-1}\cap H=1$.  The
group $G$ is \emph{CSA} if every maximal abelian subgroup is
malnormal.

\begin{rem}
If $G$ is CSA then $G$ is commutative transitive.  For, let $g\in G$
with centralizer $Z(g)$.  Consider maximal abelian $A\subset Z(g)$
and $h\in Z(g)$. Then $g\in hAh^{-1}\cap A$, so $h\in A$.  Therefore
$Z(g)=A$.
\end{rem}

\begin{lem}\label{LGs are CSA}
Limit groups are CSA.
\end{lem}
\begin{pf}
Let $H\subset G$ be a maximal abelian subgroup, consider $g\in G$,
and suppose there exists non-trivial $h\in gHg^{-1}\cap H$.  Let
$f:G\rightarrow\F$ be injective on the set
$$
\{1,g,h,[g,h]\}.
$$
Then $f([h,ghg^{-1}])=1$, which implies that $f(h)$ and
$f(ghg^{-1})$ lie in the same cyclic subgroup.  But in a free
group, this is only possible if $f(g)$ also lies in that cyclic
subgroup; so $f([g,h])=1$, and hence $[g,h]=1$.  By lemma
\ref{Abelian centralizers} it follows that $g$ commutes with every
element of $H$, so $g\in H$.
\end{pf}

\section{Embeddings in real algebraic groups}

In this section we provide solutions to exercises 8 and 9 of
\cite{BF03}, which show how to embed limit groups in real algebraic
groups and also $PSL_2(\R)$, and furthermore give some control over
the nature of the embeddings. First, we need a little real algebraic
geometry.

By, for example, proposition 3.3.13 of \cite{BCR}, every real
algebraic variety $V$ has an open dense subset
$V_\mathrm{reg}\subset V$ with finitely many connected components,
so that every component $V'\subset V_\mathrm{reg}$ is a manifold.

\begin{lem}\label{Regular components}
Consider a countable collection $V_1,V_2,\ldots\subset V$ of
closed subvarieties.  Then for any component $V'$ of
$V_\mathrm{reg}$ as above, either there exists $k$ so that
$V'\subset V_k$, or
$$
V'\cap\bigcup_{i=1}^\infty V_i
$$
has empty interior.
\end{lem}
\begin{pf}
Suppose $V'\cap\bigcup_i V_i$ doesn't have empty interior.  Then, by
Baire's Category Theorem, there exists $k$ such that $V'\cap V_k$
doesn't have empty interior.  Consider $x$ in the closure of the
interior of $V'\cap V_k$ and let $f$ be an algebraic function on
$V'$ that vanishes on $V_k$.  Then $f$ has zero Taylor expansion at
$x$, so $f$ vanishes on an open neighbourhood of $x$.  In
particular, $x$ lies in the interior of $V'\cap V_k$.  So the
interior of $V'\cap V_k$ is both open and closed, and $V'\subset
V_k$ since $V'$ is connected.
\end{pf}

\begin{lem}[Exercise 8 of \cite{BF03}]\label{Embedding in PSL2}
Let $\curlyG$ be an algebraic group over $\R$ in which $\F$ embeds.  Then for any limit
group $G$ there exists an embedding
$$
G\hookrightarrow\curlyG.
$$
In particular, $G$ embeds into $SL_2(\R)$ and $SO(3)$.
\end{lem}
\begin{pf}
Consider the variety $V=\Hom(G,\curlyG)$.  (If $G$ is of rank $r$, then $V$
is a subvariety of $\curlyG^r$, cut out by the relations of $G$.  By Hilbert's
Basis Theorem, finitely many relations suffice.)  For each
$g\in G$, consider the subvariety
$$
V_g=\{f\in V|f(g)=1\}.
$$
If $G$ does not embed into $\curlyG$ then $V$ is covered by the
subvarieties $V_g$ for $g\neq 1$.  By lemma \ref{Regular
components} every component of $V_\mathrm{reg}$ is contained in
some $V_g$, so
$$
V=V_{g_1}\cup\ldots\cup V_{g_n}
$$
for some non-trivial $g_1,\ldots,g_n\in G$.

So every homomorphism from $G$ to $\curlyG$ kills one of the
$g_i$.  But $\F$ embeds in $\curlyG$, so this contradicts the assumption
that $G$ is $\omega$-residually free.
\end{pf}

\begin{rem}
Given a limit group $G$ and an embedding $f:G\hookrightarrow
SL_2(\R)$ we have a natural map $G\to PSL_2(\R)$.  This is also an
embedding since any element in its kernel satisfies $f(g)^2=1$ and
$G$ is torsion-free.
\end{rem}

We can gain more control over embeddings into $PSL_2(\R)$ by
considering the trace function.

\begin{lem}[Exercise 9 of \cite{BF03}]
If $G$ is a limit group and $g_1,\ldots,g_n\in G$ are non-trivial
then there is an embedding $G\hookrightarrow PSL_2(\R)$ whose image
has no non-trivial parabolic elements, and so that the images of
$g_1,\ldots,g_n$ are all hyperbolic.
\end{lem}
\begin{pf}
We abusively identify each element of the variety
$V=\Hom(G,SL_2(\R))$ with the corresponding element of
$\Hom(G,PSL_2(\R))$, and call it elliptic, hyperbolic or parabolic
accordingly.  For each $g\in G$, consider the closed subvariety
$U_g$ of homomorphisms that map $g$ to a parabolic, and the open set
$W_g$ of homomorphisms that map $g$ to a hyperbolic. (Note that
$\gamma\in SL_2(\R)$ is parabolic if $|\mathrm{tr}\gamma|=2$ and
hyperbolic if $|\mathrm{tr}\gamma|<2$.) Fix an embedding
$\F\hookrightarrow SL_2(\R)$ whose image in $PSL_2(\R)$ is the
fundamental group of a sphere with open discs removed; such a
subgroup is called a \emph{Schottky group}, and every non-trivial
element is hyperbolic. Call a component $V'$ of $V_\mathrm{reg}$
\emph{essential} if its closure contains a homomorphism $G\to
SL_2(\R)$ that factors through $\F\hookrightarrow SL_2(\R)$ and maps
the $g_i$ non-trivially.

Suppose every essential component $V'$ of $V_\mathrm{reg}$ is
contained in some $U_{g'}$ for non-trivial $g'$.  Then, since
there are only finitely many components, for certain non-trivial
$g'_1,\ldots,g'_m$, every homomorphism $G\to\F$ kills one of the
$g_i$ or one of the $g'_j$, contradicting the assumption that $G$
is $\omega$-residually free. Therefore, by lemma \ref{Regular
components}, there exists an essential component $V'$ so that
$V'\cap\bigcup_{g\neq 1}U_g$ has empty interior.  In particular,
$$
\bigcap_i W_{g_i}\smallsetminus \bigcup_{g\neq 1} U_g
$$
is non-empty, as required.
\end{pf}

\section{GADs for limit groups}

In this section we provide solutions to exercises 10 and ll, and
also the related exercise 17.  For the definitions of the modular
group $\Mod(G)$, and of generalized Dehn twists, see definitions 1.6
and  1.17 respectively in \cite{BF03}.

\begin{lem}[Exercise 10 of \cite{BF03}]\label{Generalized Dehn twists}
$\Mod(G)$ is generated by inner automorphisms and generalized Dehn twists.
\end{lem}
\begin{pf}
Since the mapping class group of a surface is generated by Dehn twists (see,
for example, \cite{CB88}), it only remains to show that unimodular automorphisms of abelian vertices are generated by generalized Dehn twists.  Given such a vertex $A$, we can write
$$
G=A*_{\bar{P}(A)} B
$$
for some subgroup $B$ of $G$.  Any unimodular automorphism of $A$ is
a generalized Dehn twist of this splitting.
\end{pf}

\begin{lem}[Exercise 11 of \cite{BF03}]\label{Elliptic abelian subgroups}
Let $M$ be a non-cyclic maximal abelian subgroup of a limit group $G$.
\begin{enumerate}
\item If $G=A*_C B$ for $C$ abelian then $M$ is conjugate into either $A$
or $B$.
\item If $G=A*_C$ with $C$ abelian then either $M$ is conjugate into $A$
or there is a stable letter $t$ so that $M$ is conjugate to $M'=\langle C,t\rangle$
and $G=A*_C M'$.
\end{enumerate}
\end{lem}
\begin{pf}
We first prove 1.  Suppose $M$ is not elliptic in the splitting
$G=A*_C B$.  (Note that we don't yet know that $M$ is finitely
generated.)  Non-cyclic abelian groups have no free splittings, so
$C$ is non-trivial.  Let $T$ be the Bass--Serre tree of the
splitting.  Either $M$ fixes an axis in $T$, or it fixes a point on
the boundary.  In the latter case, there is an increasing chain of
edge groups
$$
C_1\subset C_2\subset\ldots M.
$$
But every $C_i$ is a conjugate of $C_1$, and since $M$ is malnormal
it follows that $C_i=C_1$.  So $M=C_i$, contradicting the assumption
that $M$ is not elliptic.

If $M$ fixes a line in $T$ then $M$ can be conjugated to $M'$ fixing
a line $L$ so that $A$ stabilizes a vertex $v$ of $L$ and $M'$ is of
the form $M'=C\oplus\Z$; $C$ fixes $L$ pointwise, and $\Z$ acts as
translations of $L$. Consider the edges of $L$ incident at $v$,
corresponding to the cosets $C$ and $aC$, for some $a\in
A\smallsetminus C$.  Since $aCa^{-1}=C$ and $C$ is non-trivial, it
follows from lemma \ref{LGs are CSA} that $a\in M'$. But $a$ is
elliptic so $a\in C$, a contradiction.

In the HNN-extension case, assuming $M$ is not elliptic in the
splitting we have as before that $M$ preserves a line in the
Bass--Serre tree $T$.  Conjugating $M$ to $M'$, we may assume $C$
fixes an edge in the preserved line $L$, so $C\subset M'$.  The
stabilizer of an adjacent edge is of the form $(ta)C(ta)^{-1}$,
where $a\in A$ and $t$ is the stable letter of the HNN-extension.
Therefore $C=(ta)C(ta)^{-1}$, so since $G$ is CSA it follows that
$M'=C\oplus \langle ta\rangle $ and
$$
G=A*_C M'
$$
as required.
\end{pf}

\begin{rem}\label{Remark on CSA assumptions}
Note that, in fact, the proof of lemma \ref{LGs are CSA} only used
that the vertex groups are CSA and that the edge group was maximal
abelian on one side.
\end{rem}

A one-edge splitting of $G$ is said to satisfy \emph{condition JSJ}
if every non-cyclic abelian group is elliptic in it.  Recall that a
limit group is \emph{generic} if it is freely indecomposable,
non-abelian and not a surface group.

\begin{lem}[Exercise 17 of \cite{BF03}]
If $G$ is a generic limit group then $\Mod(G)$ is generated by inner automorphisms
and generalized Dehn twists in one-edge splittings satisfying JSJ.

Furthermore, the only generalized Dehn twists that are not Dehn twists can
be taken to be with respect to a splitting of the from $G=A*_C B$ with $A=C\oplus\Z$.
\end{lem}
\begin{pf}
By lemma \ref{Generalized Dehn twists}, $\Mod(G)$ is generated by inner automorphisms
and generalized Dehn twists.  By lemma \ref{Elliptic abelian subgroups},
any splitting of $G$ as an amalgamated product satisfies JSJ.  Consider,
therefore, the splitting
$$
G=A*_C.
$$
A (generalized) Dehn twist $\delta_z$ in this splitting fixes $A$ and maps the stable letter $t\mapsto tz$, for some $z\in Z_G(C)$.

Suppose that this HNN-extension doesn't satisfy JSJ, so there exists
some (without loss, maximal) abelian subgroup $M$ that is not
elliptic in the splitting.  By lemma \ref{Elliptic abelian
subgroups}, after conjugating $M$ to $M'$, we have that
$$
G=A*_C M'
$$
and $M'=C\oplus\langle t\rangle$ where $t$ is the stable letter.  Since $M'=Z_G(C)$,
a Dehn twist $\delta_z$ along $z$ (for $c\in C$ and $n\in\Z$) fixes
$A$ and maps
$$
t\mapsto z+t.
$$
But this is a generalized Dehn twist in the amalgamated product.  So $\Mod(G)$
is, indeed, generated by generalized Dehn twists in one-edge splittings satisfying
JSJ.

Any generalized Dehn twist $\delta$ that is not a Dehn twist is in a splitting of
the form
$$
G=A*_C B
$$
with $A$ abelian, and acts as a unimodular automorphism on $A$ that preserves
$\bar{P}(A)$.  Recall that $A/\bar{P}(A)$ is finitely generated, by remark
1.15 of \cite{BF03}.  To show that it is enough to use splittings in which
$A=C\oplus\Z$, we work by induction on the rank of $A/\bar{P}(A)$.  Write
$A=A'\oplus\Z$, where $\bar{P}(A)\subset A'$.  Then there is a modular automorphism
$\alpha$, agreeing with $\delta$ on $A'$, generated by generalized Dehn twists
of the required form, by induction.  Now $\delta$ and $\alpha$ differ by
a generalized Dehn twist in the splitting
$$
G=A*_{A'} (A'*_C B)
$$
which is of the required form.
\end{pf}

\section{Constructible Limit Groups}

For the definition of CLGs see \cite{BF03}.  The definition lends
itself to the technique of proving results by a nested induction,
first on level and then on the number of edges of the GAD
$\Delta$. To prove that CLGs have a certain property, this
technique often reduces the proof to the cases where $G$ has a
one-edge splitting over groups for which the property can be
assumed.  In this section we provide solutions to exercises 12, 13, 14 and 15,
which give the first properties of CLGs, culminating in the result that all
CLGs are $\omega$-residually free.

\subsection{CLGs are CSA}

We have seen that limit groups are CSA.  This section is devoted to
proving that CLGs are also CSA.  Knowing this will prove extremely
useful in deducing the other properties of CLGs.  Note that the
property of being CSA passes to subgroups.

\begin{lem}\label{CLGs are CSA}
CLGs are CSA.
\end{lem}
By induction on the number of edges in the graph of groups $\Delta$,
it suffices to consider a CLG $G$ such that
$$
G=A*_C B
$$
or
$$
G=A*_C
$$
where each vertex group is assumed to be CSA and the edge group is
taken to be maximal abelian on one side.  (In the first case, we
will always assume that $C$ is maximal abelian in $A$.)  First, we
have an analogue of lemma \ref{Elliptic abelian subgroups}.

\begin{lem}\label{CLGs Elliptic abelian subgroups}
Let $G$ be as above.  Then $G$ decomposes as an amalgamated product
or HNN-extension in such a way that all non-cyclic maximal abelian
subgroups are conjugate into a vertex group.  Furthermore, in the
HNN-extension case we have that $C\cap C^t=1$ where $C$ is the edge
group and $t$ is the stable letter.
\end{lem}
\begin{pf}
By induction, we can assume that the vertex groups are CSA.  Note
that the proof of the first assertion of lemma \ref{Elliptic abelian
subgroups} only relies on the facts that the vertex groups are CSA
and the edge group is maximal abelian in one vertex group.  So the
amalgamated product case follows.

Now consider the case of an HNN-extension.  Suppose that, for some
stable letter $t$, $C\cap C^t$ is non-trivial.  Then since the
vertex group $A$ is commutative transitive, it follows that
$C\subset C^t$ (or $C^t\subset C$, in which case replace $t$ by
$t^{-1}$).  If $\rho:G\to G'$ is the retraction to a lower level
then, since $G'$ can be taken to be CSA, $\rho(C^t)=\rho(C)$.  But
$\rho$ is injective on edge groups so $C=C^t$ and, furthermore, $t$
commutes with $C$.

Otherwise, for every choice of stable letter $t$, $C\cap C^t$ is
trivial.  In either case, the result now follows as in the proof of
the second assertion of lemma \ref{Elliptic abelian subgroups}.
\end{pf}

Recall that a simplicial $G$-tree is \emph{$k$-acylindrical} if the
fixed point set of every non-trivial element of $G$ has diameter at
most $k$.

\begin{lem}\label{Acylindricality}
In the graph-of-groups decomposition given by lemma \ref{CLGs
Elliptic abelian subgroups}, the Bass--Serre tree is 2-acylindrical.
\end{lem}
\begin{pf}
In the amalgamated product case this is because, for any $a\in
A\smallsetminus C$, $C\cap C^a=1$.  Likewise, in the HNN-extension
case this is because $C\cap C^t=1$ where $t$ is the stable letter.
\end{pf}

\begin{pfof}{lemma \ref{CLGs are CSA}}
Let $M\subset G$ be a maximal abelian subgroup and suppose
$$
1\neq m\subset M^g\cap M.
$$
Let $T$ be the Bass--Serre tree of the splitting.  If $M$ is cyclic
then it might act as translations on a line $L$ in $T$.  Then $g$
also maps $L$ to itself.  But it follows from acylindricality that
any element that preserves $L$ lies in $M$; so $g\in M$ as required.

We can therefore assume that $M$ fixes a vertex $v$ of $T$.  If $g$
also fixes $v$ then, since the vertex stabilizers are CSA, $g\in M$.
Consider the case when
$$
G=A*_C B;
$$
the case of an HNN-extension is similar.  Then without loss of
generality $M\subset A$ and $g=ba$ for some $a\in A$ and $b\in B$ so
$M^g$ fixes a vertex stabilized by $A^b$ for some $b\in B$. Since
$C$ is maximal abelian in $A$ we have $C=M$ and since $B$ is CSA and
$m\in C\cap C^b$ it follows that $b$ commutes with $C$ so $b\in M$,
hence $g\in M$.
\end{pfof}

\subsection{Abelian subgroups}

It is no surprise that CLGs share the most elementary property of
limit groups.

\begin{lem}\label{CLGs are torsion free}
CLGs are torsion-free.
\end{lem}
\begin{pf}
The freely decomposable case is immediate by induction.  Therefore
assume $G$ is a freely indecomposable CLG of level $n$, with
$\Delta$ and $\rho:G\rightarrow G'$ as in the definition.  Suppose
$g\in G$ is of finite order.  Then $g$ acts elliptically on the
Bass--Serre tree of $\Delta$, so $g$ lies in a vertex group. Clearly
if the vertex is QH then $g=1$, and by induction if the vertex is
rigid then $g=1$. Suppose therefore that the vertex is abelian. Then
$g$ lies in the peripheral subgroup. But $\rho$ is assumed to inject
on the peripheral subgroup, so by induction $g$ is trivial.
\end{pf}

However, it is far from obvious that limit groups have abelian
subgroups of bounded rank; indeed, it is not obvious that all
abelian subgroups of limit groups are finitely generated.  But
this is true of CLGs.

\begin{lem} [Exercise 13 in \cite{BF03}] \label{Abelian subgroups of CLGs}
Abelian subgroups of CLGs are free, and there is a uniform
(finite) bound on their rank.
\end{lem}
\begin{pf}
The proof starts by induction on the level of $G$.  Let $G$ be a
CLG.  Since non-cyclic abelian subgroups have no free splittings we
can assume $G$ is freely indecomposable.  Let $\Delta$ be a
generalized abelian decomposition.  Let $T$ be the Bass--Serre tree
of $\Delta$. If $A$ fixes a vertex of $T$ then the result follows by
induction on level. Otherwise, $A$ fixes a line $T_A$ in $T$, on
which it acts by translations.  The quotient $\Delta'=T_A/A$ is
topologically a circle; after some collapses $\Delta'$ is an
HNN-extension; so the rank of $A$ is bounded by the maximum rank of
abelian subgroups of the vertex groups plus 1. So by induction the
rank of $A$ is uniformly bounded. That $A$ is free follows from
lemma \ref{CLGs are torsion free}.
\end{pf}

\subsection{Heredity}

Let $\Sigma$ be a (not necessarily compact) surface with boundary.
Then a boundary component $\delta$ is a circle or a line, and
defines up to conjugacy a cyclic subgroup
$\pi_1(\delta)\subset\pi_1(\Sigma)$.  These are called the
\emph{peripheral subgroups} of $\Sigma$.

\begin{rem} \label{Free splittings of non-compact surfaces}
Let $\Sigma$ be a non-compact surface with non-abelian fundamental
group. Then there exists a non-trivial free splitting of
$\pi_1(\Sigma)$, with respect to which all peripheral subgroups are
elliptic.
\end{rem}

\begin{lem} [Exercise 12 in \cite{BF03}] \label{Subgroups of CLGs}
Let $G$ be a CLG of level $n$ and $H$ a finitely generated subgroup.
Then $H$ is a free product of finitely many CLGs of level at most
$n$.
\end{lem}
\begin{pf}
The subgroup $H$ can be assumed to be freely indecomposable by
Grushko's theorem, so we can also assume that $G$ is freely
indecomposable.

Let $\Delta$ and $\rho:G\rightarrow G'$ be as in the definition.
Then subgroup $H$ inherits a graph-of-groups decomposition from
$\Delta$, namely the quotient of the Bass--Serre tree $T$ by $H$.
Since $H$ is finitely generated, it is the fundamental group of some
finite core $\Delta'\subset T/H$. Every vertex of $\Delta'$ covers a
vertex of $\Delta$, from which it inherits its designation as QH,
abelian or rigid.

The edge groups of $\Delta'$ are subgroups of the edge groups of
$\Delta$, so they are abelian and $\rho$ is injective on them.
Furthermore, it follows from lemma \ref{CLGs are CSA} that $H$ is
commutative transitive, so each edge group of $\Delta'$ is maximal
abelian on one side of the associated one-edge splitting.

Let $V'$ be a vertex group of $\Delta'$, a subgroup of the vertex
group $V$ of $\Delta$. There are three case to consider.
\begin{enumerate}
\item $V'\subset V$ are abelian.  Since every map $f':V'\to\Z$ with $f'(P(V'))=0$
extends to a map $f:V\to\Z$ with $f(P(V))=0$ we have that
$$\bar{P}(V')\subset\bar{P}(V)$$ so $\rho$ is injective on
$\bar{P}(V')$.
\item $V'\subset V$ are QH.  If $V'$ is of infinite index in $V$ then $V'$ is
the fundamental group of a non-compact surface, so by remark
\ref{Free splittings of non-compact surfaces} $H$ is freely
decomposable. Therefore it can be assumed that $V'$ is of finite
degree $m$ in $V$.  In particular, $V$ is the fundamental group of a
compact surface that admits a pseudo-Anosov automorphism.
Furthermore, let $g,h\in V$ be such that $\rho([g,h])\neq 1$. Then
because CLGs are commutative transitive,
$$
\rho([g^m,h^m])\neq 1.
$$
But $g^m,h^m\in V'$, so $\rho(V')$ is non-abelian.
\item $V'\subset V$ are rigid.  Then $\tilde{V}'\subset
\tilde{V}$ because CLGs are commutative transitive, so $\rho|_{V'}$ is injective.
\end{enumerate}
Therefore $\rho|_H:H\rightarrow G'$ and $\Delta'$ satisfy the
properties for $H$ to be a CLG.
\end{pf}

\subsection{Coherence}

A group is \emph{coherent} if every finitely generated subgroup is
finitely presented.  Note that free groups and free abelian groups
are coherent.   For limit groups, coherence is an instance of a more
general phenomenon, as in the next lemma.  Recall that a group is
\emph{slender} if every subgroup is finitely generated. Finitely
generated abelian groups are slender.

\begin{lem} \label{Graphs of coherent groups}
The fundamental group of a graph of groups with coherent vertex
groups and slender edge groups is coherent.
\end{lem}
\begin{pf}
Let $\Delta$ be a graph of groups, with coherent vertex groups and
slender edge groups.  Let $G=\pi_1(\Delta)$ and $H\subset G$ a
finitely generated subgroup.  Then $H$ inherits a graph-of-groups
decomposition from $\Delta$ given by taking the quotient of the
Bass--Serre tree $T$ of $\Delta$ by the action of $H$.  Since $H$ is
finitely generated it is the fundamental group of some finite core
$\Delta'\subset T/H$.  But, by induction on the number of edges,
$H=\pi_1(\Delta')$ is finitely presented.
\end{pf}

\begin{lem}[Exercise 12 in \cite{BF03}] \label{CLGs are coherent}
CLGs are coherent, in particular finitely presented.
\end{lem}
\begin{pf}
In the case of a free decomposition the result is immediate.  In the
other case, the (free abelian) edge groups of $\Delta$ are finitely
generated, so slender and coherent, by lemma \ref{Abelian subgroups
of CLGs}. Therefore all vertex groups are finitely generated; in
particular, abelian vertex groups are coherent.  Finitely generated
surface groups are also coherent. Rigid vertex groups embed into a
CLG of lower level, so by lemma \ref{Subgroups of CLGs} they are
free products of coherent groups and hence coherent by induction.
The result now follows by lemma \ref{Graphs of coherent groups}.
\end{pf}

\subsection{Finite $K(G,1)$}

That CLGs have finite $K(G,1)$ follows from the fact that graphs
of aspherical spaces are aspherical.

\begin{thm}[Proposition 3.6 of \cite{SW}] \label{Aspherical graphs of groups}
Let $\Delta$ be a graph of groups; suppose that for every vertex
group $V$ there exists finite $K(V,1)$, and for every edge group $E$
there exists a finite $K(E,1)$.   Then for $G=\pi_1(\Delta)$, there
exists a finite $K(G,1)$.
\end{thm}

Surface groups and abelian groups have finite Eilenberg--Mac~Lane
spaces.  Rigid vertices embed into a CLG of lower level, so by lemma
\ref{Subgroups of CLGs} and induction they also have finite
Eilenberg--Mac~Lane spaces.

\begin{cor} [Exercise 13 in \cite{BF03}] \label{CLGs have finite K(G,1)}
If $G$ is a CLG then there exists a finite $K(G,1)$.
\end{cor}

\subsection{Principal cyclic splittings}

A \emph{principal cyclic splitting} of $G$ is a one-edge splitting
of $G$ with cyclic edge group, such that the image of the edge group
is maximal abelian in one of the vertex groups; further, if it is an
HNN-extension then the edge group is required to be maximal abelian
in the whole group.  The key observation about principal cyclic
splittings is that any non-cyclic abelian subgroup is elliptic with
respect to them---in other words, they are precisely those cyclic
splittings that feature in the conclusion of lemma \ref{CLGs
Elliptic abelian subgroups}.  Applying lemma \ref{CLGs Elliptic
abelian subgroups}, to prove that every freely indecomposable,
non-abelian CLG has a principal cyclic splitting it will therefore
suffice to produce any non-trivial cyclic splitting (since we now
know that CLGs are CSA).

\begin{prop}[Exercise 14 in \cite{BF03}] \label{Principal cyclic splittings of CLGs}
Every non-abelian, freely indecomposable CLG admits a principal
cyclic splitting.
\end{prop}
\begin{pf}
Let $G$ be a CLG.  As usual, by induction it suffices to consider
the cases when $G$ splits as an amalgamated product or
HNN-extension.  It suffices to exhibit any cyclic splitting of $G$,
as observed above.

Suppose
$$
G=A*_C B.
$$
If $C$ is cyclic the result is immediate, so assume $C$ is
non-cyclic abelian.  If either vertex group is freely decomposable
then so is $G$, since $C$ has no free splittings; if both vertex
groups are abelian then so is $G$.  Therefore $A$, say, is freely
indecomposable and non-abelian so has a principal cyclic splitting,
which we shall take to be of the form
$$
A=A'*_{C'} B'.
$$
(It might also be an HNN-extension, but this doesn't affect the
proof.)  Because it is principal $C$ is conjugate into a vertex, say
$B'$; so $G$ now decomposes as
$$
G=A'*_{C'} (B'*_C B).
$$
which is a cyclic splitting as required.

The proof when $G=A*_C$ is the same.
\end{pf}

\subsection{A criterion in free groups}

To prove that a group $G$ is $\omega$-residually free, it suffices
to show that for any finite $X\subset G\smallsetminus 1$ there
exists a homomorphism $f:G\rightarrow\F$ with $1\notin f(X)$.  So a
criterion to show that an element of $\F$ is not the identity will
be useful.

\begin{lem}\label{Free group criterion}
Let $z\in\F\smallsetminus 1$, and consider an element $g$ of the
form
$$
g=a_0z^{i_1} a_1z^{i_2}a_2\ldots a_{n-1}z^{i_n}a_n
$$
where $n\geq 1$ and, whenever $0<k<n$, $[a_k,z]\neq 1$. Then $g\neq
1$ whenever the $|i_k|$ are sufficiently large.
\end{lem}

Choose a generating set for $\F$ so the corresponding Cayley graph
is a tree $T$.  An element $u\in\F$ specifies a geodesic
$[1,u]\subset T$. Likewise, a string of elements
$u_0,u_1,\ldots,u_n\in\F$ defines a path
$$
[1,u_0]\cdot u_0[1,u_1]\cdot\ldots\cdot (u_0\ldots u_{n-1})[1,u_n]
$$
in $T$, where $\cdot$ denotes concatenation of paths.  The key
observation we will use is as follows.  The length of a word
$w\in\F$ is denoted by $|w|$.

\begin{rem}\label{Key observation}
Suppose $z$ is cyclically reduced and has no proper roots.  Let
$a\in\F$ be such that $a$ and $az$ both lie in $L\subset T$ the axis
of $z$. If $j$ is minimal such that $z^j$ lies in the geodesic
$[a,az]$ then, setting $u=z^ja^{-1}$ and $v=az^{1-j}$, it follows
that $uv=z=vu$; in particular, either $u$ or $v$ is trivial and
$[a,z]=1$.
\end{rem}

\begin{pfof}{lemma \ref{Free group criterion}}
It can be assumed that $z$ is cyclically reduced and has no proper
roots.

Assume that, for each $k$, $|z^{i_k}|\geq |a_{k-1}|+|a_k|+|z|$. Let
$L\subset T$ be the axis of $z$. Denote by $g_k$ the partial product
$$
g_k=a_0z^{i_1} a_1z^{i_2}a_2\ldots a_{k-1}z^{i_k}a_k.
$$
The path $\gamma$ corresponding to $g$ is of the form
$$
[1,a_0]\cdot g_0[1,z^{i_1}]\cdot g_0z^{i_1}[1,a_1]\cdot\ldots\cdot
g_{n-1}[1,z^{i_n}]\cdot g_{n-1}z^{i_n}[1,a_n].
$$
Suppose that $g=1$ so this path is a loop.  Each section of the form
$g_k[1,z^{i_k}]$ lies in a translate of $L$, the axis of $z$.  Since
$T$ is a tree, for at least one such section $\gamma$ enters and
leaves $g_kL$ at the same point---otherwise $\gamma$ is a
non-trivial loop.  Since $|z^{i_k}|>|a_{k-1}|+|a_k|+|z|$ it follows
that both $g_kz^{i_k}a_k$ and $g_kz^{i_k}a_kz$ lie in $g_kL$ and so
$[a_k,z]=1$ by remark \ref{Key observation}.
\end{pfof}

\subsection{CLGs are limit groups}

\begin{thm}[Exercise 15 in \cite{BF03}] \label{CLGs are LGs}
CLGs are $\omega$-residually free.
\end{thm}

Since the freely decomposable case is immediate, let $\Delta$, $G'$
and $\rho$ be as in the definition of a CLG in \cite{BF03}. By
induction, $G'$ can be assumed $\omega$-residually free.  As a warm
up, and for use in the subsequent induction, we first prove the
result in the case of abelian and surface vertices.

\begin{lem} \label{Abelian case}
Let $A$ be a free abelian group and $\rho:A\rightarrow G'$ a
homomorphism to a limit group. Suppose $P\subset A$ is a nontrivial
subgroup of finite corank closed under taking roots, on which $\rho$
is injective. Then for any finite subset $X\subset A\smallsetminus1$
there exists an automorphism $\alpha$ of $A$, fixing $P$, so that
$1\notin \rho\circ\alpha(X)$.
\end{lem}
\begin{pf}
Since $\ker\rho$ is a subgroup of $A$ of positive codimension, for
given $x\in A\smallsetminus 0$ a generic automorphism $\alpha$
certainly satisfies $\alpha(x)\notin\ker\rho$.   Since $X$ is
finite, therefore, there exists $\alpha$ such that
$\alpha(x)\notin\ker\rho$ for any $x\in X$.
\end{pf}

We now consider the surface vertex case.

\begin{prop}\label{Surface case}
Let $S$ be the fundamental group of a surface $\Sigma$ with
non-empty boundary, with $\chi(\Sigma)\leq -1$, and
$\rho:S\rightarrow G'$ a homomorphism injective on each peripheral
subgroup and with non-abelian image. Then for any finite subset
$X\subset S\smallsetminus1$ there exists an automorphism $\alpha$ of
$S$, induced by an automorphism of $\Sigma$ fixing the boundary
components pointwise, such that $1\notin\rho\circ\alpha(X)$.
\end{prop}

Let the surface $\Sigma$ have $b>0$ boundary components and Euler
characteristic $\chi<0$.  When $\Sigma$ is cut along a two-sided
simple closed curve $\gamma$, the resulting pieces either have lower
genus (defined to be $1-\frac{1}{2}(\chi+b)$) or fundamental groups
of strictly lower rank, depending on whether $\gamma$ was separating
or not.  The simplest cases all have fundamental groups that are
free of rank 2.

\begin{exa}[The simplest cases]\label{Simplest cases}
Suppose $S$ is free of rank 2.  By lemma \ref{2-generator LGs},
$\rho(S)$ is free or free abelian; but $\rho(S)$ is assumed
non-abelian, so $\rho(S)$ is free and $\rho$ is injective.
\end{exa}

For the more complicated cases, the idea is to find a suitable
simple closed curve $\zeta$ along which to cut to make the surface
simpler. In order to apply the proposition inductively, $\zeta$
needs the following properties:
\begin{enumerate}
\item $\rho(\zeta)\neq 1$;
\item the fundamental group $S'$ of any component of $\Sigma\smallsetminus\zeta$
must have $\rho(S')$ non-abelian.
\end{enumerate}

Let's find this curve in some examples.

\begin{figure}[htp]
\begin{center}
\psfrag{1}{$d_1$} \psfrag{2}{$d_2$} \psfrag{3}{$d_3$}
\psfrag{4}{$d_4$}\psfrag{5}{$d_1d_2$}
\includegraphics[width=0.7\textwidth]{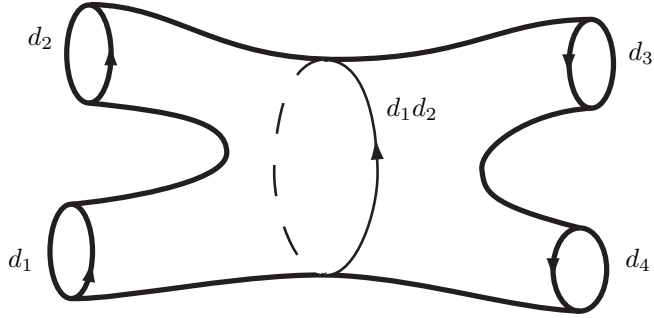}
\caption{A four-times punctured sphere}
\end{center}
\end{figure}

\begin{exa}[Punctured spheres] \label{Punctured spheres}
Suppose $\Sigma$ is a punctured sphere, so
$$
S=\langle d_1,\ldots,d_n|\prod_j d_j\rangle.
$$
Assume $n\geq 4$ and $\rho(d_i)\neq 1$ for all $i$. Define a
relation on $\{1,\ldots,n\}$ by
$$
i\sim j \Leftrightarrow \rho([d_i,d_j])= 1.
$$
Since $G'$ is commutative transitive, $\sim$ is an equivalence
relation.  Because the image is non-abelian, there are at least
two equivalence classes.  Since
$$
\prod_i d_i=1
$$
any equivalence class has at least two elements in its complement.
Relabelling if necessary, it can now be assumed that
$$
\rho([d_1,d_2]),\rho([d_3,d_4])\neq 1.
$$
Now if the boundary curves have been coherently oriented then $d_1
d_2$ has a representative that is a simple closed curve. Take
$\zeta$ as this representative.
\end{exa}

The case when $\Sigma$ is non-orientable is closely related.

\begin{exa}[Non-orientable surfaces]\label{Non-orientable surfaces}
Suppose $\Sigma$ is non-orientable so
$$
S=\langle c_1,\ldots,c_m,d_1,\ldots,d_n|\prod_i
c_i^2\prod_jd_j\rangle.
$$
Exactly the same argument as in the case of a punctured sphere would
work if it could be guaranteed that $\rho(c_i)\neq 1$ for all $i$.

Fix some $c_k$, therefore, and suppose $\rho(c_k)=1$.  Let $\gamma$
be a simple closed curve representing $c_k$.  Then $d_1c_k$ has a
representative $\delta$ which is a simple closed curve, and
$\rho(d_1c_k)\neq 1$.  Furthermore, $\Sigma\smallsetminus\gamma$ and
$\Sigma\smallsetminus\delta$ are homeomorphic surfaces, and a
homeomorphism between them extends to an automorphism of $\Sigma$
mapping $\gamma$ to $\delta$.   This homeomorphism can be chosen not
to alter any of the other $c_i$ or the $d_j$.

Therefore, after an automorphism of $\Sigma$, it can be assumed
that $\rho(c_i)\neq 1$ for all $i$, so a suitable $\zeta$ can be
found as in the previous example.
\end{exa}

\begin{figure}[htp]
\begin{center}
\psfrag{1}{$d_1$} \psfrag{2}{$d_2$} \psfrag{3}{$d_3$}
\psfrag{a}{$a_1$}\psfrag{b}{$b_1$}\psfrag{a1}{$a_1d_1$}
\includegraphics[width=0.7\textwidth]{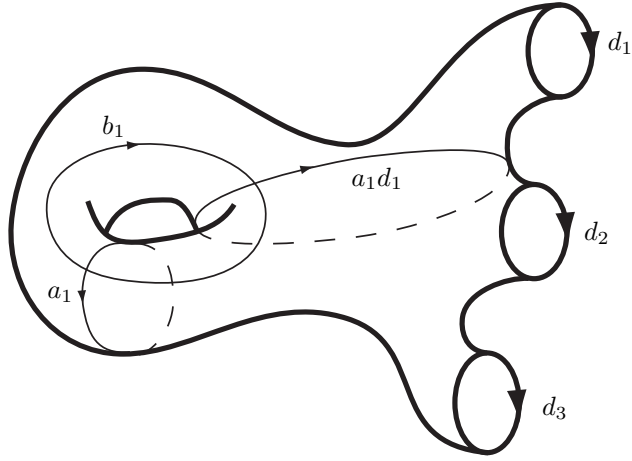}
\caption{The positive-genus case.}
\end{center}
\end{figure}

\begin{exa}[Positive-genus surfaces]\label{Positive genus}
Suppose $\Sigma$ is an orientable surface of positive genus, so
$$
S=\langle
a_1,b_1,\ldots,a_g,b_g,d_1,\ldots,d_n|\prod_i[a_i,b_i]\prod_j d_j
\rangle.
$$
Assume that $g,n\geq 1$.  If, for example, $\rho(a_1)\neq 1$ then
$\zeta$ can be taken to be a simple closed curve representing $a_1$.
Otherwise, $\rho(a_1d_1)\neq 1$ and $a_1d_1$ has a simple closed
representative.  It remains to show that the single component of
$\Sigma\smallsetminus\zeta$ has non-abelian image.

Cutting along $\zeta$ expresses $S$ as an HNN-extension:
$$
S=S'*_\Z.
$$
Let $t$ be the stable letter, and suppose $\rho(S')$ is abelian.
Then $\zeta\in tS't^{-1}\cap S'$ so in particular
$\rho(tS't^{-1}\cap S')$ is non-trivial.  But $G'$ is a limit group
and hence CSA, so $\rho(t)$ commutes with $\rho(S')$, contradicting
the assumption that $\rho(S)$ is non-abelian.
\end{exa}

Note that examples \ref{Punctured spheres}, \ref{Non-orientable
surfaces} and \ref{Positive genus} cover all the more complicated
surfaces with boundary.
\begin{pfof}{proposition \ref{Surface case}}
Example \ref{Simplest cases} covers all the simplest cases. Suppose
therefore that $\Sigma$ is more complicated.   To apply the
inductive hypothesis, an essential simple closed curve
$\zeta\in\Sigma$ is needed such that $\rho(\zeta)\neq 1$ and, for
any component $S'$ of $\Sigma\smallsetminus\zeta$, $\rho(S')$ is
non-abelian. This is provided by examples \ref{Punctured spheres},
\ref{Non-orientable surfaces} and \ref{Positive genus}.

For simplicity, assume $\zeta$ is separating.  The non-separating
case is similar. Then $\Sigma\smallsetminus\zeta$ has two
components, $\Sigma_1$ and $\Sigma_2$.  Let $S_i=\pi_1(\Sigma_i)$,
and denote by $X_i$ the \emph{syllables} of $X$ in $S_i$---the
elements of $S_i$ that occur in the normal form of some $x\in X$
with respect to the splitting over $\langle\zeta\rangle$. Because
the pieces $\Sigma_i$ are simpler than $\Sigma$ there exists
$\alpha\in\Aut_0(\Sigma)$ and $f:G'\rightarrow\F$ such that
$$
1\notin f\circ\rho\circ\alpha([\zeta,X_1\cup X_2]).
$$
Consider $\xi\in X$.  The proposition follows from the claim that,
for all sufficiently large $k$,
$$
f\circ\rho\circ\delta_\zeta^k\circ\alpha(\xi)\neq 1
$$
where $\delta_\zeta$ is a Dehn twist in $\zeta$.  If $\xi$ is a
power of $\zeta$ then the result is immediate. Otherwise, with
respect to the one-edge splitting of $G$ over $\langle\zeta\rangle$,
$\xi$ has reduced form
$$
\sigma_0\tau_0\sigma_1\tau_1\ldots\sigma_n\tau_n
$$
where the $\sigma_i\in X_1$ and the $\tau_j\in X_2$.  The image
$x^{(k)}=f\circ\rho\circ\delta_\zeta^k\circ\alpha(\xi)$ is of the
form
$$
z^ks_0z^{-k}t_0z^ks_1z^{-k}t_1\ldots z^ks_nz^{-k}t_n
$$
where $z=f\circ\rho(\zeta)$, $s_i=f\circ\rho\circ\alpha(\sigma_i)$
and $t_i=f\circ\rho\circ\alpha(\tau_i)$.  This expression for
$x^{(k)}$ satisfies the hypotheses of lemma \ref{Free group
criterion}, so $x^{(k)}\neq 1$ for all sufficiently large $k$.
\end{pfof}

The proof of theorem \ref{CLGs are LGs} is very similar to the proof
of proposition \ref{Surface case}. The theorem follows from the
following proposition, by induction on level.

\begin{prop}\label{Heart of the matter}
Let $G$ be a freely indecomposable CLG, let $G'$ be
$\omega$-residually free, and let $\Delta$ and  $\rho$ be as usual.
For any finite subset $X\subset G\smallsetminus1$ there exists a
modular automorphism $\alpha$ of $G$ such that $1\notin
\rho\circ\alpha(X)$.
\end{prop}
\begin{pf}
As usual, the proposition is proved by induction on the number of
edges of $\Delta$.  The case of $\Delta$ having no edges follows
from lemmas \ref{Abelian case} and \ref{Surface case}, and the fact
that $\rho$ is injective on rigid vertices.  By induction on level,
$G'$ is a limit group.

Now suppose $\Delta$ has an edge group $E$.  For simplicity, assume
$E$ is separating.  The non-separating case is similar. Then
removing the edge corresponding to $E$ divides $\Delta$ into two
subgraphs $\Delta_1$ and $\Delta_2$.  Let $G_i=\pi_1(\Delta_i)$, and
denote by $X_i$ the syllables of $X$ in $G_i$.  Without loss assume
$E$ is maximal abelian in $G_1$. Fix non-trivial $\zeta\in E$. By
induction there exists $\alpha\in\Mod(\Delta)$ and
$f:G'\rightarrow\F$ such that
$$
1\notin f\circ\rho\circ\alpha([\zeta,X_1]\cup X_2).
$$
Consider $\xi\in X$.  The proposition follows from the claim that,
for all sufficiently large $k$,
$$
f\circ\rho\circ\delta_\zeta^k\circ\alpha(\xi)\neq 1.
$$
If $\xi\in E$ then the result is immediate.  Otherwise, with
respect to the one-edge splitting of $G$ over $E$, $\xi$ has
reduced form
$$
\sigma_0\tau_0\sigma_1\tau_1\ldots\sigma_n\tau_n
$$
where the $\sigma_i\in X_1$ and the $\tau_j\in X_2$.  The image
$x^{(k)}=f\circ\rho\circ\delta_\zeta^k\circ\alpha(\xi)$ is of the
form
$$
z^ks_0z^{-k}t_0z^ks_1z^{-k}t_1\ldots z^ks_nz^{-k}t_n
$$
where $z=f\circ\rho(\zeta)$, $s_i=f\circ\rho\circ\alpha(\sigma_i)$
and $t_i=f\circ\rho\circ\alpha(\tau_i)$.  In particular, canceling
across those $t_i$ that commute with $z$, we have
$$
x^{(k)}=u_0z^{k\epsilon_1}u_1\ldots u_{n-1}z^{k\epsilon_n}u_n
$$
where $\epsilon_i=\pm1$ and $u_i$ don't commute with $z$ for
$0<i<n$. This second expression for $x^{(k)}$ satisfies the
hypotheses of lemma \ref{Free group criterion}, so $x^{(k)}\neq 1$
for all sufficiently large $k$.
\end{pf}

\section{The Shortening Argument}\label{Shortening argument section}

We consider a sequence of $G$-trees $T_i$, arising from
homomorphisms $f_i:G\to\F$, that converge in the Gromov topology to
a $G$-tree $T$.  By the results of section 3 of \cite{BF03}, if the
action of $G$ on the limit tree $T$ is faithful then it gives rise
to a generalized abelian decomposition for $G$.  This section is
entirely devoted to the solution of exercise 16, which is
essentially Rips and Sela's shortening argument---an ingenious means
of using this generalized abelian decomposition to force the action
on the limit tree to be unfaithful.

\subsection{Preliminary ideas}\label{Shortening argument subsection}

Once again, fix a generating set for $\F$ so that the corresponding
Cayley graph is a tree, and let $|w|$ denote the length of a word
$w\in\F$. Fix a generating set $S$ for $G$.  For $f:G\rightarrow\F$,
let
$$
|f|=\max_{g\in S}|f(g)|.
$$
A homomorphism is \emph{short} if
$$
|f|\leq |\iota\circ f\circ\alpha|
$$
whenever $\alpha$ is a modular automorphism of $G$ and $\iota$ is
an inner automorphism of $\F$.

\begin{thm}[Exercise 16 of \cite{BF03}]\label{Shortening argument}
Suppose every $f_i$ is short.  Then the action on $T$ is not
faithful.
\end{thm}
The proof is by contradiction.  We assume therefore, for the rest of
section \ref{Shortening argument section}, that the action is
faithful.  By the results summarized in section 3 of \cite{BF03},
the action of $G$ on $T$ gives a GAD $\Delta$ for $G$. The idea is,
if $T_i$ are the limiting trees with basepoints $x_i$, to construct
modular automorphisms $\phi_n$ so that
$$
d_i(x_i,f_i\circ\phi_i(g)x_i)<d_i(x_i,f_i(g)x_i)
$$
for all sufficiently large $i$.  Then apply these automorphisms to
carefully chosen basepoints.

All constructions of the limit tree $T$, such as the asymptotic cone
\cite{vdDW}, use some form of based convergence: basepoints $x_i\in
T_i$ are fixed, and converge to a basepoint $[x_i]\in T$.  Because
the $f_i$ are short,
$$
\max_{g\in S}d_\F(1,f_i(g))\leq\max_{g\in S}d_\F(t,f_i(g)t)
$$
for all $t\in T_\F$; otherwise, conjugation by the element of $\F$
nearest to $t$ leads to a shorter equivalent homomorphism. It
follows that $1\in T_i$ is always a valid basepoint; we set
$x=[1]\in T$ to be the basepoint for $T$.

The proof of theorem \ref{Shortening argument} goes on a
case-by-case basis, depending on whether $[x,gx]$ intersects a
simplicial part or a minimal part of $T$.

\subsection{The abelian part}

The next proposition is a prototypical shortening result for a
minimal vertex.

\begin{prop}\label{Proposition 5.2}
Let $V$ be an abelian vertex group of $\Delta$.  For $g\in G$, let
$l(g)$ be the translation length of $g$ on $T$. Fix $\epsilon>0$.
Then for any finite subset $S\subset V$ there exists a modular
automorphism $\phi$ of $G$ such that
$$
\max_{g\in S}l(\phi(g))<\epsilon.
$$
\end{prop}
\begin{pf}
The minimal $V$-invariant subtree $T_V$ is a line in $T$, on which
$V$ acts indiscretely.  Since $S$ is finite, $V$ can be assumed
finitely generated.  It suffices to prove the theorem in the case
where $S$ is a basis for $V$.  Assume furthermore that each element
of $S$ translates $T_V$ in the same direction.

Suppose the action of $V$ on $T_V$ is free.  Let
$S=\{g_1,\ldots,g_n\}$, ordered so that
$$
l(g_1)>l(g_2)>\ldots>l(g_n)>0.
$$
Since the action is indiscrete, there exists an integer $\lambda$
such that
$$
l(g_1)-\lambda l(g_2)<\frac{1}{2}l(g_2).
$$
Applying the automorphism that maps $g_1\mapsto g_1-\lambda g_2$ and
proceeding inductively, we can make $l(g_1)$ as short as we like.

If the action of $V$ is not free then $V=V'\oplus V_0$ where $V'$
acts freely on $T_V$ and $V_0$ fixes $T_V$ pointwise.  Applying the
free case to $V'$ gives the result.
\end{pf}

The aim is to prove the following theorem.

\begin{thm}\label{Abelian shortening}
Let $V$ be an abelian vertex. Then for any finite subset $S\subset
G$ there exists a modular automorphism $\phi$ such that for any
$g\in S$:
\begin{enumerate}
\item if $[x,gx]$ intersects a translate of $T_V$ in a segment of positive length
then
$$
d(x,\phi(g)x)<d(x,gx);
$$
\item otherwise, $\phi(g)=g$.
\end{enumerate}
\end{thm}
\begin{pf}
By a result of J. Morgan (claim 3.3 of \cite{Mo}---the article is
phrased in terms of laminations), the path $[x,gx]$ intersects
finitely many translates of $T_V$ in non-trivial segments.  Let
$\epsilon$ be the minimal length of all such segments across all
$g\in S$.  Assume that $g\in S$ is such that $[x,gx]$ intersects a
translate of $T_V$ non-trivially.

Suppose first that $x$ lies in a translate of $T_V$, so without loss
of generality $x\in T_V$.  Then $g$ has a non-trivial decomposition
in the GAD provided by corollary 3.16 of \cite{BF03} of the form
$$
g=a_0b_1a_1\ldots a_n
$$
where the $a_i$ lie in $V$ and the $b_i$ are products of elements of
other vertices and loop elements.  Write $g_i=a_0b_1\ldots
b_{i-1}a_{i-1}$.  The decomposition can be chosen so that each
component of the geodesic $[x,gx]$ that lies in $g_iT_V$ is
non-trivial.  For each $i$, decompose $[x,b_ix]$ as
$$
[x,s_i]\cdot[s_i,t_i]\cdot[t_i,b_ix]
$$
where $[x,s_i]$ and $[t_i,b_ix]$ are maximal segments in $T_V$ and
$b_iT_V$ respectively.  Then
$$
[x,gx]=[x,a_0s_1]\cdot g_1[s_1,t_1]\cdot
g_1[t_1,b_1a_1s_2]\cdot\ldots\cdot g_n[s_n,t_n]\cdot
g_n[t_n,b_na_nx]\\
$$
where each $[s_i,t_i]$ and $[t_i,b_ia_is_{i+1}]$ is a non-trivial
segment.  Therefore
\begin{eqnarray*}
d(x,gx) & = & d(x,a_0s_1)+\sum_{i=1}^n
d(s_i,t_i)+\sum_{i=1}^{n-1}d(t_i,b_ia_is_{i+1})+d(t_n,b_na_nx)\\
& \geq & \sum_{i=1}^n d(s_i,t_i) + (n+1)\epsilon.
\end{eqnarray*}

Since $V$ acts indiscretely on the line $T_V$, by modifying the
$b_i$ by elements of $V$ it can be assumed that
$$
d(x,s_i),d(t_i,b_ix)<\frac{1}{4}\epsilon.
$$
By proposition \ref{Proposition 5.2} there exists $\phi\in\Mod(G)$
such that $\phi(b_i)=b_i$ for all $b_i$ and
$$
d(x,\phi(a_i)x)<\frac{1}{2}\epsilon
$$
for all $a_i$.  Now as before $[x,\phi(g)x]$ decomposes as
$$
[x,\phi(a_0)s_1]\cdot
\phi(g_1)[s_1,t_1]\cdot\ldots\cdot\phi(g_n)[t_n,b_n\phi(a_n)x]
$$
so
\begin{eqnarray*}
d(x,\phi(g)x) & = & d(x,\phi(a_0)s_1)+\sum_{i=1}^n
d(s_i,t_i)+\sum_{i=1}^{n-1}d(t_i,b_i\phi(a_i)s_{i+1})\\& &+d(t_n,b_n\phi(a_n)x)\\
& < & d(x,\phi(a_0)x)+\sum_{i=1}^n
d(s_i,t_i)+\sum_{i=1}^{n-1}d(b_ix,b_i\phi(a_i)x)\\& &+d(b_nx,b_n\phi(a_n)x)+\frac{n}{2}\epsilon\\
& < & \sum_{i=1}^n d(s_i,t_i)+ \big(n+\frac{1}{2}\big)\epsilon.
\end{eqnarray*}
Therefore, $d(x,\phi(g)x)<d(x,gx)-\frac{1}{2}\epsilon$ and in
particular the result follows.

Now suppose $x$ does not lie in a translate of $T_V$.  Then $g$ has
a non-trivial decomposition in the corresponding GAD of the form
$$
g=b_0a_1b_1\ldots b_n
$$
where the $a_i$ lie in $V$ and the $b_i$ are products of elements of
other vertices and loop elements.  Let $g'=a_1b_1\ldots b_{n-1}a_n$.
Let $x'$ be the first point on $[x,gx]$ in a translate of $T_V$, so
$x'=b_0y\in b_0T_V$ for some $y\in T_V$.  Likewise let $x''$ be the
last point on $[x,gx]$, so $x''=b_0g'z\in b_0g'T_V$ for some $z\in
T_V$.  Since the action of $V$ on $T_V$ is indiscrete we can modify
$b_n$ by an element of $V$ and assume that
$d(y,z)<\frac{1}{4}\epsilon$.

Then the geodesic $[x,gx]$ decomposes as
$$
[x,gx]=[x,b_0y]\cdot[b_0y,b_0g'z]\cdot[b_0g'z,gx]
$$
so
$$
d(x,gx)> d(x,b_0y)+d(y,g'y)+d(z,b_nx)-\frac{1}{4}\epsilon.
$$
Applying the first case to $g'$ and $y$ we obtain $\phi\in\Mod(G)$
such that
$$
d(y,\phi(g')y)<d(y,g'y)-\frac{1}{2}\epsilon
$$
so
\begin{eqnarray*}
d(x,\phi(g)x)&<&d(x,b_0y)+d(y,\phi(g')y)+d(z,b_nx)+\frac{1}{4}\epsilon\\
&<&
d(x,b_0y)+d(y,g'y)-\frac{1}{2}\epsilon+d(z,b_nx)+\frac{1}{4}\epsilon\\
&<& d(x,gx)
\end{eqnarray*}
as required.
\end{pf}

\subsection{The surface part}

The surface part is dealt with by Rips and Sela, in \cite{RS94}, in
the following theorem.

\begin{thm}[Theorem 5.1 of \cite{RS94}]\label{Surface shortening}
Let $V$ be a surface vertex. Then for any finite subset $S\subset G$
there exists a modular automorphism $\phi$ such that for any $g\in
S$:
\begin{enumerate}
\item if $[x,gx]$ intersects a translate $T_V$ in a segment of positive length
then
$$
d(x,\phi(g)x)<d(x,gx);
$$
\item otherwise, $\phi(g)=g$.
\end{enumerate}
\end{thm}

Rips and Sela use the notion of groups of interval exchange
transformations, which are equivalent to surface groups, and prove
an analogous result to proposition \ref{Proposition 5.2}.  The rest
of the proof is the same as that of theorem \ref{Abelian
shortening}.

\subsection{The simplicial part}

It remains to consider the case where $[x,gx]$ is contained in the
simplicial part of $T$.

\begin{thm}\label{Simplicial shortening}
Let $S\subset G$ be finite and let $x\in T$.  Then there exist
$\phi_n\in\Mod(\Delta)$ such that, for all $g\in S$,
$$
d(x,\phi_n(g)x)=d(x,gx);
$$
furthermore, for all $g\in S$ that do not fix $x$, and for all
sufficiently large $n$,
$$
d_n(x_n,f_n\circ\phi_n(g)x_n)<d_n(x_n,f_n(g)x_n).
$$
\end{thm}

Let $e$ be a closed simplicial edge containing $x$.  The proof of
the theorem is divided into cases, depending on whether the image of
$e$ is separating in $T/G$.  In both cases, the following lemma will
prove useful.

\begin{lem}\label{Defining the mn}
Let $A$ be a vertex group of the splitting over $e$.   Let $T_A$ be
the minimal $A$-invariant subtree of $T$; conjugating $A$, we can
assume that $T_A\cap e$ is precisely one point, $y$.  Fix any
non-trivial $c\in C=\Stab(e)$. Then there exists a sequence of
integers $m_n$ such that, for any $a\in A$,
$$
d_n(x_n,f_n(c^{-m_n}ac^{m_n})x_n)\rightarrow d(y,ay)
$$
as $n\rightarrow \infty$.
\end{lem}
\begin{pf}
The key observation is that
$$
2d_n(x_n,\Axis(f_n(c)))<d_n(x_n,f_n(c)x_n)\rightarrow 0,
$$
and the same holds for the $y_n$.  Let $x'_n$ be the nearest point
on $\Axis(f_n(c))$ to $x_n$; likewise, let $y'_n$ be the nearest
point on $\Axis(f_n(c))$ to $y_n$.  Then, for each $n$, there exists
$m_n$ such that
$$
d_n(f_n(c^{m_n})x'_n,y'_n)<l(f_n(c))\rightarrow 0
$$
as $n\rightarrow\infty$.  Therefore
$$
d_n(f_n(c^{m_n})x_n,y_n)\rightarrow 0
$$
as $n\rightarrow\infty$, and the result follows.
\end{pf}

The next lemma helps with the case where the image of $e$ is
separating.

\begin{lem}\label{Separating automorphism}
Assume the image of $e$ is separating, so the induced splitting is
$$
G=A*_CB.
$$
Assume furthermore that, with the notation of the previous lemma,
$x\neq y$.  Then there exists $\alpha_n\in\Mod(\Delta)$ such that,
for all $g\in S$:
\begin{enumerate}
\item if $g\in A$ then $\alpha_n(g)=g$;
\item if $g\notin A$ then
$$
d_n(x_n,f_n\circ\alpha_n(g)x_n)<d_n(x_n,f_n(g)x_n).
$$
\end{enumerate}
\end{lem}
\begin{pf}
Fix a non-trivial $c\in C$, and let $\delta_c$ be the Dehn twist in
$c$ that is the identity when restricted to $A$.  Let
$\alpha_n=\delta_c^{m_n}$ where $m_n$ are the integers given by
lemma \ref{Defining the mn}. Any $g\notin A$ has normal form
$$
g=a_0b_1a_1\ldots b_la_l
$$
with the $a_i\in A\smallsetminus C$ and the $b_i\in B\smallsetminus
C$, except for $a_0$ and $a_l$ which may be trivial.  Therefore
$d(x,gx)=\sum_i d(x,b_ix)+\sum_i d(x,a_ix)$.  Fix $\epsilon>0$.  If
$a_i$ is non-trivial then $a_i\notin C$ and so
$$
d(x,a_ix)=2d(x,y)+d(y,a_iy).
$$
Let $k$ be the number of $a_i$ that are non-trivial (so $l-1\leq
k\leq l+1$).  Therefore, for all sufficiently large $n$,
$$
d_n(x_n,f_n(g)x_n)>\sum_i d(x,b_ix) + \sum_i
d(y,a_iy)+2kd(x,y)-\epsilon.
$$
By contrast, for all sufficiently large $n$,
$$
d_n(x_n,f_n\circ\alpha_n(g)x_n)<\sum_i d(x,b_ix)+\sum_i
d(y,a_iy)+\epsilon
$$
by lemma \ref{Defining the mn}.  By assumption $x\neq y$, so
$d(x,y)>0$. Therefore taking $\epsilon<kd(x,y)$ gives the result.
\end{pf}

We now turn to the non-separating case.

\begin{lem}\label{Defining the pn}
Assume the image of $e$ is non-separating, so the splitting induced
by $e$ is
$$
G=A*_C.
$$
Let $t$ be a stable letter.  As before, conjugate $A$ so that
$T_A\cap e$ is precisely one point $y$.   Fix any non-trivial $c\in
C=\Stab(e)$. Then there exists a sequence of integers $p_n$ such
that
$$
d_n(y_n,f_n(tc^{p_n})y_n)\rightarrow 0
$$
as $n\rightarrow \infty$.  Therefore, for any fixed integer $j$,
$$
d_n(y_n,f_n(tc^{p_n})^j y_n)\rightarrow 0
$$
as $n\rightarrow \infty$.
\end{lem}
\begin{pf}
As in the proof of lemma \ref{Defining the mn}, by the definition of
Gromov convergence,
$$
2d_n(y_n,\Axis(f_n(c)))<d_n(y_n,f_n(c)y_n)\rightarrow 0
$$
as $n\rightarrow \infty$, and similarly,
$$
2d_n(f_n(t^{-1})y_n,\Axis(f_n(c)))\rightarrow 0
$$
as $n\rightarrow \infty$.  Let $y'_n$ be the nearest point on
$\Axis(f_n(c))$ to $y_n$, and let $y''_n$ be the nearest point on
$\Axis(f_n(c))$ to $f_n(t^{-1})y_n$.  Then there exist integers
$p_n$ such that
$$
d_n(f_n(c^{p_n})y'_n,y''_n)\rightarrow 0
$$
as $n\rightarrow\infty$.  The result now follows.
\end{pf}

\begin{lem}\label{Non-separating automorphism}
Assume the situation is in lemma \ref{Defining the pn}.  Then there
exists $\alpha_n\in\Mod(\Delta)$ such that, for all $g\in S$:
\begin{enumerate}
\item if $g\in C$ then $\alpha_n(g)=g$;
\item if $g\notin C$ then
$$
d_n(x_n,f_n\circ\alpha_n(g)x_n)<d_n(x_n,f_n(g)x_n).
$$
\end{enumerate}
\end{lem}
\begin{pf}
Fix a stable letter $t$ that translates $x$ away from $y$.  Fix a
non-trivial $c\in C$ and let $i_c\in\Mod(G)$ be conjugation by $c$.
Set $\alpha_n=i_c^{m_n}\circ\delta_c^{p_n}$, where $m_n$ are
integers given by lemma \ref{Defining the mn} and $p_n$ are given by
lemma \ref{Defining the pn}. Any $g$ is of the form
$$
g=a_0t^{j_1}a_1\ldots t^{j_l}a_l
$$
with $j_i\neq 0$ and the $a_i\in A\smallsetminus C$ except for $a_0$
and $a_l$ which may be trivial.  Unlike in the case of a separating
edge, we have to be a little more careful in estimating $d(x,gx)$
because the natural path from $x$ to $gx$ given by the decomposition
of $g$ may backtrack.  To be precise, backtracking occurs when
$a_i\neq 1$ and $j_{i+1}<0$ and also when $j_i>0$ and $a_{i+1}\neq
1$.  Let $k$ be the number of $i$ for which backtracking does not
occur, so $0\leq k\leq 2$. Then
$$
d(x,gx)=\sum_i d(y,a_iy)+\sum_i d(x,t^{j_i}x)+2kd(x,y).
$$
Fix $\epsilon>0$. Then for all sufficiently large $n$,
$$
d_n(x_n,f_n(g)x_n)>\sum_i d(y,a_iy)+2kd(x,y)+\sum_i
d(x,t^{j_i}x)-\epsilon.
$$
Now for each $i$,
$$
d_n(f_n(c^{m_n})x_n,f_n((tc^{p_n})^{j_i})f_n(c^{m_n})x_n)\rightarrow
0
$$
and
$$
d_n(f_n(c^{m_n})x_n,f_n(a_ic^{m_n})x_n)\rightarrow d(y,a_iy).
$$
So for all sufficiently large $n$,
$$
d_n(x_n,f_n\circ\alpha_n(g)x_n)<\sum_i d(y,a_iy)+\epsilon.
$$
Taking $2\epsilon<2kd(x,y)+\sum_i d(x,t^{j_i}x)$ gives the result.
\end{pf}

We are now ready to prove the theorem.

\begin{pfof}{theorem \ref{Simplicial shortening}} Suppose first that $x$ lies in the interior of
an edge $e$.  If $e$ has separating image in the quotient then lemma
\ref{Separating automorphism} can be applied both ways round, giving
rise to modular automorphisms $\alpha_n$ and $\beta_n$.  The theorem
is then proved by taking $\phi_n=\alpha_n\circ\beta_n$.  If $e$ is
non-separating then applying lemma \ref{Non-separating automorphism}
and taking $\phi_n=\alpha_n$ gives the result.

Suppose now that $x$ is a vertex.  For each orbit of edges $[e]$
adjoining $x$, let $\alpha_n^e$ be the result of applying lemma
\ref{Separating automorphism} or lemma \ref{Non-separating
automorphism} as appropriate to $e$. Now taking
$$
\phi_n=\alpha_n^{e_1}\circ\ldots\circ\alpha_n^{e_p}
$$
where $[e_1],\ldots,[e_p]$ are the orbits adjoining $x$ gives the
required automorphism.\end{pfof}

This is the final piece of the shortening argument.

\begin{pfof}{theorem \ref{Shortening argument}}
Fix a generating set $S$ for $G$.  Let $f_i:G\to\F$ be a sequence of
short homomorphisms corresponding to the convergent sequence of
$G$-trees $T_i$.  Let $T$ be the limiting $G$-tree and suppose that
the action of $G$ on $T$ is faithful.  By corollary 3.16 of
\cite{BF03} this induces a GAD $\Delta$ for $G$.  Let $x\in T$ be
the basepoint fixed in subsection \ref{Shortening argument
subsection}.

Composing the automorphisms given by theorems \ref{Abelian
shortening} and \ref{Surface shortening} there exists
$\alpha\in\Mod(\Delta)$ such that, for any $g\in G$,
$$
d(x,\phi(g)x)<d(x,gx)
$$
if $[x,gx]$ intersects an abelian or surface component of $T$ and
$\phi(g)=g$ otherwise.  By theorem \ref{Simplicial shortening}, for
all sufficiently large $i$ there exist $\beta_i\in\Mod(\Delta)$ such
that $d(x,\beta_i(g)x)=d(x,gx)$ and, furthermore,
$$
d_i(1,f_i\circ\beta_i(g))<d_i(1,f_i(g))
$$
whenever $[x,gx]$ is a non-trivial arc in the simplicial part of the
tree.  It follows that for $\phi_i=\beta_i\circ\alpha$,
$$
d_i(1,f_i\circ\phi_i(g))<d_i(1,f_i(g))
$$
for all $g\in S$ and all sufficiently large $i$.  This contradicts
the assumption that the $f_i$ were short.
\end{pfof}

\section{Bestvina and Feighn's geometric approach}

In section 7 of \cite{BF03}, Bestvina and Feighn provide a more
geometric proof of their Main Proposition.  In this section we
provide proofs of the exercises needed in this argument.

\subsection{The space of laminations}

Recall that $\ML(K)$ is the space of measured laminations on $K$,
and $\Proj\ML(K)$ is its quotient by the action of $\R_+$.  Let
$E$ be the set of edges of $K$.

\begin{prop}[Exercise 18 of \cite{BF03}]
The space of measured laminations on $K$ can be identified with a
closed cone in $\R_+^E-\{0\}$, given by the triangle inequality
for each 2-cell of $K$.  Hence, when $\ML(K)$ is endowed with the
corresponding topology, $\Proj\ML(K)$ is compact.
\end{prop}
\begin{pf}
Recall that two laminations are considered equivalent if they
assign the same measure to each edge.  Therefore it suffices to
show existence of a lamination with the prescribed values on the
edges.  First, for each edge $e$ with $\int_e \mu>0$, fix a closed
proper subinterval $I_e$ contained in the interior of $e$.  Now
fix a Cantor function $c_e:I_e\rightarrow [0,\int_e\mu]$.  This
gives a measure $\mu$ on $e$, given by
$$
\int_J\mu=\int_{I_e\cap J} c_e d\lambda
$$
where $d\lambda$ is Lebesgue measure on $\R$.  Now suppose $e_1,
e_2,e_3$ are the edges of a simplex in $K$.  Divide $e_1$ into
intervals $e_1^2$ and $e_1^3$ so that
$$
2\int_{e_1^2} d\mu=\int_{e_1}d\mu +\int_{e_2}d\mu - \int_{e_3}d\mu
$$
and $e_1^2$ shares a vertex with $e_2$, and similarly for $e_1^3$.
Divide $e_2$ and $e_3$ likewise. Fix a Cantor set in each $e_i^j$.
Now for each distinct $i,j$ inscribe a lamination between $e_i^j$
and $e_j^i$.  Since any path transverse to this lamination can be
homotoped to an edge path respecting the lamination, the measure on
the edges determines a transverse measure to the lamination.
\end{pf}

\subsection{Matching resolutions in the limit}

A measured lamination on $K$ defines a $G$-tree.  The next exercise
shows the close relation between the topology on the space of
laminations and the topology on the space of trees.  For the
definition of a resolution, see \cite{BF03}.  The solution is most
easily phrased in terms of some explicit construction of the
limiting tree.  I shall use the asymptotic cone, $T_\omega$; $T$ can
be realized as the minimal $G$-invariant subtree of $T_\omega$.  For
the definition of the asymptotic cone see, for example, \cite{vdDW}.
To see how to choose basepoints and scaling to ensure that the
action is non-trivial see, for example, \cite{P88}.

\begin{prop}[Exercise 19 in \cite{BF03}]
Consider $f_i$-equivariant resolutions
$$
\phi_i:\tilde{K}\to T_\F.
$$
Suppose $\lim T_{f_i}=T$, $\lim \Lambda_{\phi_i}=\Lambda$ and the
sequences $(|f_i|)$ and $(||\phi_i||_K)$ are comparable.  Then
there is a resolution that sends lifts of leaves of $\Lambda$ to
points of $T$ and is a Cantor function on edges of $\tilde{K}$.
\end{prop}
\begin{pf}
A resolution $\phi:\tilde{K}\rightarrow T$ is determined by a
choice of $\phi(\tilde{v})$ for a lift $\tilde{v}$ of each vertex
$v$ of $K$.

First, define a resolution $\phi':\tilde{K}\to T_\omega$ by setting
$\phi'(\tilde{v})=[\phi_i(\tilde{v})]$.  Since $(|f_i|)$ and
$(||\phi_i||_K)$ are comparable, $\phi'(\tilde{v})$ is a valid point
of $T_\omega$.  The resolution $\phi'$ maps leaves of $\Lambda$ to
points, and is a Cantor function on edges.  However, $T_\omega$ is
far from minimal. Let $\pi:T_\omega\to T$ be closest-point
projection to the minimal invariant subtree, which is equivariantly
isomorphic to $T$.  Now let $\phi=\pi\circ\phi'$'; this is a
resolution that still maps leaves of $\Lambda$ to points, and is a
Cantor function on edges, as required.
\end{pf}

\subsection{Finding kernel elements carried by leaves}

Exercise 20 of \cite{BF03} relies heavily on the results of
\cite{BF95}. The most important result is a structure theorem for
resolutions of stable actions on real trees, summarized in the
following theorem.

\begin{thm}[Theorems 9.4 and 9.5 of \cite{BF95}]\label{FP actions on trees}
Let $\Lambda$ be a lamination on a 2-complex $K$, resolving a
stable action of $G=\pi_1(K)$ on a real tree $T$.  Then
$$
\Lambda=\Lambda_1\sqcup\ldots\sqcup\Lambda_k.
$$
Each component has a standard neighbourhood $N_i$ carrying a
subgroup $H_i\subset G$.  Let $T_i$ be the minimal $H_i$-invariant
subtree of $T$.  Each component is of one of the following types.
\begin{enumerate}
\item \textbf{Surface type}.  $N_i$ is a cone-type 2-orbifold, with some
annuli attached.  $H_i$ fits into a short exact sequence
$$
1\rightarrow \ker T_i\rightarrow H_i\rightarrow\pi_1(O)\rightarrow
1
$$
where $O$ is a cone-type 2-orbifold.
\item \textbf{Toral type}.  $T_i$ is a line, and $H_i$ fits into a short
exact sequence
$$
1\rightarrow \ker T_i\rightarrow H_i\rightarrow A\rightarrow 1
$$
where $A\subset\Isom(\R)$.
\item \textbf{Thin type}.  $H_i$ splits over an arc stabilizer, carried by
a leaf of $\Lambda_i$.
\item \textbf{Simplicial type}.  All the leaves of $\Lambda_i$ are compact,
and $N_i$ is an interval bundle over a leaf.  $H_i$ fits into a
short exact sequence satisfying
$$
1\rightarrow\ker T_i\rightarrow H_i\rightarrow C\rightarrow 1
$$
where $C$ is finite.
\end{enumerate}
Furthermore, if $E$ is a subgroup carried by a leaf, $E$ fits into a
short exact sequence of the form
$$
1\rightarrow\kappa\rightarrow E\rightarrow C\rightarrow 1
$$
where $\kappa$ fixes an arc of $T$ and $C$ is finite or cyclic.
\end{thm}

In particular, the standard neighbourhoods induce a
graph-of-spaces decomposition for $K$, and a corresponding
graph-of-groups decomposition for $G$.  The vertex spaces are the
$N_i$ and the closures of the components of $K-\cup_iN_i$. The
edge spaces are boundary components of the $N_i$, and are all
contained in a leaf. See theorem 5.13 of \cite{BF95}.

The proof of this exercise will also make use of the following
result.

\begin{prop}[Corollary 5.9 of \cite{BF95}]\label{Stability}
If $h\in H_i$ fixes an arc of $T_i$ then $h\in\ker T_i$.
\end{prop}

We are now ready to prove the exercise.

\begin{thm}[Exercise 20 of \cite{BF03}]
In the situation of the exercise, the lamination $\Lambda$ has a
leaf carrying non-trivial elements of the kernel.
\end{thm}
\begin{pf}
Note that $G/\ker T$ is a limit group.  Suppose no elements of $\ker
T$ are carried by a leaf of $\Lambda$.

Consider $\Gamma$ the graph of groups for $G$ induced by
$\Lambda$.   The aim is to show that $\Gamma$ really is a GAD.
Since a GAD decomposition can be used to shorten, this contradicts
the assumption that the $f_i$ are short.  We deal with each sort
of vertex in turn.

\begin{enumerate}
\item Suppose $\Lambda_i$ is of surface type.  Then $N_i$ is a
cone-type 2-orbifold, with some annuli attached. Suppose $g\in H_i$
is carried by an annulus.  Then $g$ fixes an arc of $T_i$, so by
proposition \ref{Stability}, $g\in\ker T_i$.  But $T_i$ contains a
tripod, and tripod stabilizers are trivial, so $g\in\ker T$
contradicting the assumption.  Therefore $N_i$ can be assumed to
have no attached annuli.  Consider an element $g\in H_i$ carried by
the leaf corresponding to a cone-point.  Then $g$ has finite order,
so $g\in\ker T$, since $G/\ker T$ is a limit group.  This
contradicts the assumption, so $N_i$ has no cone-points.  Therefore
$N_i$ is genuinely a surface.  Moreover, $N_i$ carries a
pseudo-Anosov homeomorphism, since it carries a minimal lamination.

\item If $\Lambda_i$ is toral, then $H_i$ is an extension
$$
1\rightarrow\ker T_i\rightarrow H_i\rightarrow A\rightarrow 1
$$
for $A\subset\Isom\R$.  The elements of $\ker T_i$ are carried by
annuli in $N_i$.  But $\ker T_i$ itself fits into an exact
sequence
$$
1\rightarrow\kappa\rightarrow\ker T_i\rightarrow A'\rightarrow 1
$$
where $\kappa\subset\ker T$ and $A'$ is abelian.  In order not to
contradict the assumption that no elements of the kernel are carried
by a leaf, therefore, $\kappa$ must be trivial; so we have
$$
1\rightarrow A'\rightarrow H_i\rightarrow A\rightarrow 1.
$$
and $H_i$ acts faithfully on $T$.  In particular, $H_i$ embeds in
the limit group $G/\ker T$, and so is a limit group.  But $A'$ is
normal; since limit groups are torsion-free and CSA, it follows
that $H_i$ is free abelian.

\item If $\Lambda_i$ is thin, then $G$ splits over a subgroup $H$
fixing an arc of $T_i$.  By proposition \ref{Stability}, $H\subset
\ker{T_i}$.  But $T_i$ contains a tripod, and tripod stabilizers
are trivial, so $H\subset\ker T$; since $H$ is carried by a leaf,
$H$ must be trivial by assumption.  But this contradicts the
assumption that $G$ is freely indecomposable.
\item If $\Lambda_i$ is simplicial, then $H_i$ fits into the short
exact sequence
$$
1\rightarrow\ker T_i\rightarrow H_i\rightarrow C\rightarrow 1.
$$
for finite $C$.  As in the toral case, the assumption implies that
$\ker T_i$ is abelian, and $H_i$ embeds in $G/\ker T$, and so is a
limit group.  But, again, $H_i$ is torsion-free and CSA; so $C$ is
trivial, and $H_i$ is abelian and fixes an arc of $T$.
\end{enumerate}

Now consider an edge-group $E$ of $\Gamma$.  Then $E$ is carried
by a leaf, and satisfies
$$
1\rightarrow\kappa\rightarrow E\rightarrow C\rightarrow 1
$$
where $C$ is cyclic and $\kappa$ fixes an arc of $T$.  Then
$\kappa$ fits into a short exact sequence
$$
1\rightarrow\kappa'\rightarrow\kappa\rightarrow A\rightarrow 1
$$
where $\kappa'\subset\ker T$ and $A$ is abelian.  By assumption,
therefore, $\kappa'$ is trivial and $\kappa$ is abelian;
furthermore, $E$ acts faithfully on $T$, so embeds in $G/\ker T$
and is a limit group.   Therefore $E$ is free abelian.

In conclusion, $\Gamma$ is a GAD.   Just as in the proof of
theorem \ref{Shortening argument}, this contradicts the assumption
that the $f_i$ are all short.
\end{pf}

\subsection{Examples of limit groups}

To complete their argument, Bestvina and Feighn need some elementary
examples of limit groups.  This theorem is required.

\begin{thm}[Exercise 21 of \cite{BF03}]
Let $\Delta$ be a 1-edged GAD of a group $G$ with a homomorphism $q$
to a limit group $\Gamma$.  Suppose:
\begin{enumerate}
\item the vertex groups of $\Delta$ are non-abelian,
\item the edge group of $\Delta$ is maximal abelian in each vertex
group, and
\item $q$ is injective on vertex groups of $\Delta$.
Then $G$ is a limit group.
\end{enumerate}
\end{thm}

This theorem is just a special case of proposition \ref{Heart of the
matter}.

\bibliographystyle{plain}
\bibliography{exercises}

\begin{thebibliography}{10}

\bibitem{BF95}
M.~Bestvina and M.~Feighn.
\newblock Stable actions of groups on real trees.
\newblock {\em Invent. Math.}, 121, 1995.

\bibitem{BF03}
M.~Bestvina and M.~Feighn.
\newblock Notes on {S}ela's work: Limit groups and {M}akanin-{R}azborov
  diagrams, 2003.
\newblock Preprint.

\bibitem{BCR}
J.~Bochnak, M.~Coste, and M-F. Roy.
\newblock {\em G\'{e}om\'{e}trie alg\'{e}brique r\'{e}elle}, volume~3 of {\em
  Ergebnisse der {M}athematik und ihrer {G}renzgebiete}.
\newblock Springer-Verlag, Berlin, 1987.

\bibitem{CB88}
Andrew~J. Casson and Steven~A. Bleiler.
\newblock {\em Automorphisms of surfaces after {N}ielsen and {T}hurston}.
\newblock Number~9 in London {M}athematical {S}ociety {S}tudent {T}exts.
  {C}ambridge {U}niversity {P}ress, 1988.

\bibitem{L59}
R.~S. Lyndon.
\newblock The equation $a^2b^2=c^2$ in free groups.
\newblock {\em Michigan Math. J}, 6, 1959.

\bibitem{Mo}
J.~Morgan.
\newblock Ergodic theory and free actions of groups on {R}-trees.
\newblock {\em Invent. Math.}, 94, 1988.

\bibitem{P88}
F.~Paulin.
\newblock Topologie de {G}romov \'equivariante, structures hyperboliques et
  arbes r\'eels.
\newblock {\em Invent. Math.}, 94(1), 1988.

\bibitem{RS94}
E.~Rips and Z.~Sela.
\newblock Structure and rigidity in hyperbolic groups {I}.
\newblock {\em GAFA}, 4(3), 1994.

\bibitem{SW}
Peter Scott and Terry Wall.
\newblock Topological methods in group theory.
\newblock In {\em Proc. Sympos., Durham, 1977}, volume~36 of {\em London Math.
  Soc. Lecture Note Ser.}, pages 137--203. Cambridge University Press, 1979.

\bibitem{Se1}
Z.~Sela.
\newblock Diophantine geometry over groups {I}: {M}akanin-{R}azborov diagrams.
\newblock {\em Publ. Inst. Hautes \'Etudes Sci.}, 93, 2001.

\bibitem{vdDW}
L.~van~den Dries and A.~Wilkie.
\newblock On {G}romov's theorem concerning groups of polynomial growth and
  elementary logic.
\newblock {\em J. Algebra}, 89, 1984.

\end{thebibliography}
\end{document}